\documentclass[a4paper,psamsfonts,12pt]{amsart}
\usepackage{amsmath,amssymb}
\usepackage[totalwidth=480pt,totalheight=680pt]{geometry}
\def\ifnextchar#1#2#3{\let\tmpnce=#1%
    \def\tmpnca{#2}\def\tmpncb{#3}\futurelet\tmpncc\ifnch}%
  \def\ifnch{\ifx\tmpncc\tmpnce\let\tmpncd=\tmpnca%
	\else\let\tmpncd=\tmpncb\fi\tmpncd}
\def\Arg{\operatorname{Arg}}

   \def\pdiff#1#2{\dfrac{\partial#1}{\partial#2}}

\def\thei{{\bold i}}

\def\modulus#1{\left|#1\right|}
\def\monthname{\ifcase\month\or Jan\or Feb\or March\or Apr\or %
    May\or June\or July\or Aug\or Sept\or Oct\or Nov\or Dec\fi}
\def\norm#1{\left\|#1\right\|}
      \def\normp#1_#2{\norm{#1}_{#2}}

\def\Set#1{\left\{\,#1\,\right\}}

\let\C\Complex

\let\R\Real

\let\sphere\Sphere

\def\map#1#2#3{#1~\colon~#2~\to~#3}
\input epsf.tex
\title[Off-center Dynamics]{The Dynamics of Off-center Reflection}
\author[Thomas Au]{Thomas Kwok-keung Au}
\address{Dept of Mathematics, The Chinese University of Hong Kong, Shatin, Hong Kong.}
\email{thomasau@cuhk.edu.hk}

\newtheorem{theorem}{Theorem}[section]
\newtheorem{proposition}[theorem]{Proposition}
\newtheorem{corollary}[theorem]{Corollary}

\newtheorem{lemma}{Lemma}

\def\tR{\tilde{R}}

\begin{document}
\begin{abstract}
Dynamical properties of a two-parameter circle map, called off-center
reflection, are studied.  Certain symmetry breaking phenomena in
the bifurcation process are illustrated and discussed.
\end{abstract}
\maketitle
\thispagestyle{empty}
\setlength\baselineskip{26pt}
\setlength\parskip{2ex}
\setlength\parindent{0em}
\pagestyle{headings}

\section*{Background and Definition}
We study the dynamics of a two-parameter family of circle maps
$
\map{R_{r,\Omega}} {\sphere^1}{\sphere^1}
$
called the off-center reflection.  When $\Omega=\pi$, this map is
a one-dimensional analog of the general map raised in \cite[problem~21]{Yau},
which is geometrically a reflection in the circle.  For other values
of $\Omega$, it can be seen as a reflection with a deviation between
the reflected and incident angles.
Iterations of this map are not the natural sucessive reflections
in the circle; nevertheless, this map is interesting for various
reasons.
The off-center reflection can also be seen as a perturbation (with
small $r>0$) of rotation by $\Omega$. 
In fact, it has an analytical form which extends the well-known Arnold
circle map, \cite{Arnold}.
Its dynamics is related to the perturbation properties of
Mathieu type differential equation, \cite{Ar2}.
Furthermore, when the perturbation parameter $r$ goes to~1, the map
goes to another famous circle map, the doubling map.

The off-center reflection is introduced in \cite{AL} by the following
geometric description.  Fix a point $L$ inside the unit circle $\sphere^1$.
For a point $\phi\in\sphere^1$, a ray is emitted from $L$ to $\phi$.
This ray is ``reflected'' to hit $\sphere^1$ again
at a point, denote $R_{r,\Omega}(\phi)$ in the future.  This point
is defined to be the image of $\phi$ under the map.
It is quoted ``reflected'' because the ``reflected'' angle has a
constant deviation $\frac{\pi-\Omega}2$ from the incident angle
$\iota(\phi)$.  To have the map geometrically well-defined, there is
some restriction on $\Omega$.  However, we will see later that it is
analytically meaningful for other $\Omega$.  In fact, it is
sufficient to consider $\Omega \in (-\pi,\pi]$.  Furthermore, since
the action has certain symmetries, it is no loss
of generality to assume the point source $L$ at $(r,0)$ with $0\leq
r < 1$.  This is why the off-center reflection is given by the two
parameters $(r,\Omega)$.
\begin{center}
\mbox{\epsfysize=49mm \epsfbox{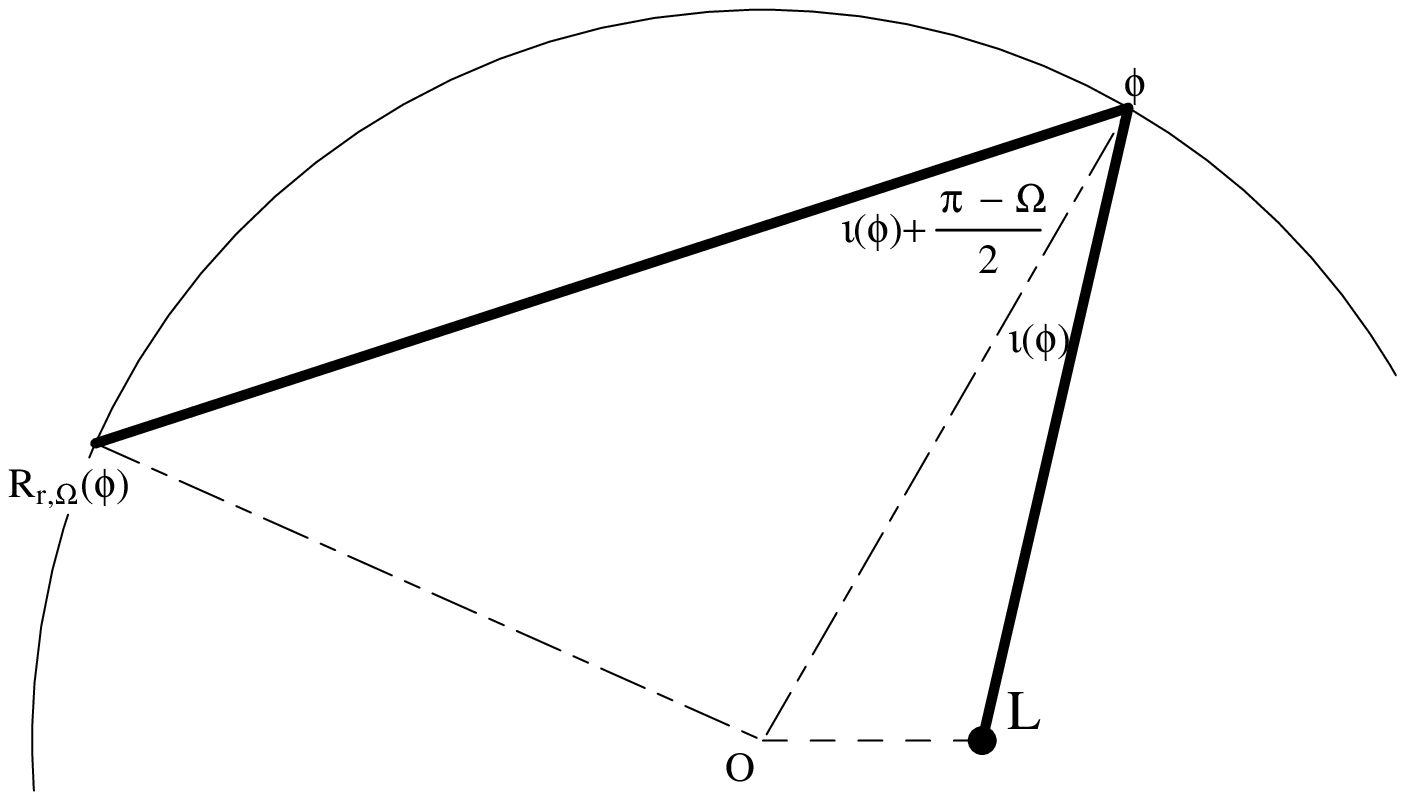}}
\end{center}

In \cite{AL}, particular interest is placed on $R_{r,\pi}$, i.e.,
when the reflected angle equals the incident angle.
There, the link between dynamics and contact geometry of the map
is studied.  Here, we deal with the dynamical properties of the family
with some focus on $\Omega = 0, \pi$.  It is because in these two
cases, the off-center reflection map repects the symmetry of the
circle and symmetric cycles may occur.  We are particularly interested
in when symmetric cycles of the map break into asymmetric ones.
The properties presented in this article may be considered the first
steps to understand the dynamics of the map.  It is expected that
deeper studies may bring forth more understanding to
general circle maps.

In \S 1, we will introduce the basic analytical properties of the map.
In \S 2, the attracting orbits, especially the symmetric ones, of the
map are investigated.  Then, finally in \S 3, the bifurcation of the
map is looked into.  In particular, we give an explanation of how and
when the symmetric orbits go through a period preserving pitch-fork
bifurcation.
The analysis in \cite{Br} of similar behavior among certain cubic 
polynomials is borrowed.

\setcounter{section}{0}
\section{Analysis}
For our purpose, we consider $\sphere^1 \subset
\C\simeq\R^2$ and it is covered by $\R$ under the exponential map
$x\mapsto e^{\thei x}$.
Then the off-center reflection
$$
\phi \mapsto R_{r,\Omega}(\phi) :
\sphere^1 \longrightarrow \sphere^1
$$
has a unique continuous lifting,
$\map{\tilde{R}_{r,\Omega}} {\R}{\R}$, which takes $0$ to $\Omega$. 
Since $\tilde{R}_{r,\Omega}(x+2\pi) = \tilde{R}_{r,\Omega}(x) +
2\pi$, we may focus our attention on the interval $(-\pi,\pi]$.
Let the incident angle be denoted by $\iota(x)$ for $x \in (-\pi,\pi]$.
Since $\iota(x) \to 0$ as $x\to \pm\pi$,
it defines a continuous $2\pi$-periodic function on $\R$, which is also
denoted as $\iota(x)$.
In fact, it can be written in terms of the principal argument (with
values between $-\pi$ and $\pi$) as
$$
\iota(x) = \Arg(\cos x - r + \thei\sin x) - x;
$$
and ${\tilde R}_{r,\Omega}(x) = x + \Omega - 2\iota(x)$.
The first few derivatives of $\tilde R$ are listed below.
\begin{align*}
\iota'(x) &= \frac {r(\cos x - r)}{(\cos x-r)^2 + \sin^2 x};
\\
\tilde{R}_{r,\Omega}'(x) &= 1 - 2\iota'(x) =
\frac {1-4r\cos x + 3r^2} {(\cos x-r)^2 + \sin^2 x} \\
\tilde{R}_{r,\Omega}''(x) &= \frac {2r(1-r^2)\sin x} {\left[(\cos x-r)^2 + \sin^2 x\right]^2} \\
\tilde{R}_{r,\Omega}'''(x) &= \frac {2r(1-r^2) \left[(1+r^2)\cos x - 2r(1+\sin^2x)\right]} {\left[(\cos x-r)^2 + \sin^2 x\right]^3}.
\end{align*}

This incident angle has a series expression given in \cite{AL},
$\displaystyle
\iota(x) = \sum_{k=1}^\infty \frac{r^k}{k} \sin(kx).
$
Therefore, the lift of $R_{r,\Omega}$ mapping $0$ to $\Omega$ is given by
\begin{align*}
\tR_{r,\Omega}(x) &= x + \Omega - 2\iota(x) \\
&= x + \Omega - 2\sum_{k=1}^\infty \frac{r^k}{k} \sin(kx).
\end{align*}
From this, we see that the Arnold circle map $x \mapsto x +
\Omega - \varepsilon\sin x$ can be seen as a truncated
version of the off-center reflection.
The technique of our analysis in this paper may also be
adapted to a similar study on the Arnold circle map.

We will first establish a nice property for the dynamics of
a map.  The dynamics of the off-center reflection in the
range $0\leq r < 1/3$ is relatively simple.  The following
proposition is useful for studying the range $r > 1/3$.
\begin{proposition}
For $1/3<r<1$ and all $\Omega$, $R_{r,\Omega}$ is a map with negative
Schwarzian derivative.
\end{proposition}
\begin{proof}
The Schwarzian derivative of $\tilde{R}_{r,\Omega}$ is given by
$$
\frac{\tR'''_{r,\Omega}(x)}{\tR'_{r,\Omega}(x)} -
\frac{3}{2} \left(
\frac{\tR''_{r,\Omega}(x)}{\tR'_{r,\Omega}(x)} \right)^2 =
\frac {r(1-r^2)H(r,\cos x)} {(1-4r\cos x+3r^2)^2 (1-2r\cos x+r^2)^2}
$$
where
$$
\begin{aligned}
H(r,\cos x) &=
(2+28r^2+6r^4)\cos x \\
&\qquad\qquad {}- r\left[ 13+19r^2-(1-r^2)\cos 2x+4r\cos 3x\right]
\\
&=
-(14r+18r^3) + (2+40r^2+6r^4)\cos x \\
&\qquad\qquad {}+
2r(1-r^2)\cos^2x - 16r^2\cos^3x \\
\end{aligned}
$$
Let $y=\cos x$ and consider $H(r,y)$ for $y \in [-1,1]$. 
We have
\begin{align*}
H(r,-1) &= -2(1+r)^3(1+3r) \\
&= -2 - 12r - 24r^2 - 20r^3 - 6r^4 < 0; \\
H(r,1) &= 2(1-r)^3(1-3r) < 0, \qquad \text{for $1/3 < r < 1$.}
\end{align*}
Moreover, its derivative wrt $y$ satisfies
\begin{align*}
\partial_y H &= (2+40r^2+6r^4) + 4r(1-r^2)y - 48r^2y^2; \\
\intertext{and for $1/3 < r < 1$,}
\partial_y H(r,-1) &= 2(1+r)(1-r^2)(1-3r) < 0, \\
\partial_y H(r,1) &= 2(1-r)(1-r^2)(1+3r) > 0.
\end{align*}
We see that, inside $(-1,1)$, the cubic polynomial $y
\mapsto H(r,y)$ can only have a single critical point and it
is a minimum.  Hence $H(r,y) \leq \max\{ H(r,-1), H(r,1) \}
< 0$.
\end{proof}

\section{Attracting Orbits}
It is easy to see from our earlier formula that $\tilde{R}_{r,\Omega}'(x)
\geq 0$ for $0 \leq r \leq \frac{1}{3}$, with equality only when $r =
\frac{1}{3}$ and $\cos x = 1$.  Thus, $R_{r,\Omega}$ is a homeomorphism
for $0\leq r \leq 1/3$, so the dynamics is trivial.  On the other hand,
$R_{r,\Omega}$ is only a degree~1 map for $r > 1/3$.  We would like
to explore the dynamics of it in the coming sections.

In the study of periodic orbits of $R_{r,\Omega}$, the information about
the function $\iota$ is often helpful.  Since $\iota$ is a $2\pi$-periodic
odd function, it is sufficient to know its properties in the interval
$[0, \pi]$.  To be precise, $\iota(x) > 0$ and is concave down for
$x\in(0,\pi)$.  Its maximum value of $\pi/2-a_r$ is attained at $a_r$
where $a_r$ is the angle satisfying $0\leq a_r \leq \pi/2$ and $\cos a_r = r$. 
A picture of its graph will be helpful to see its properties.
\begin{center}
\mbox{\epsfysize=4cm \epsfbox{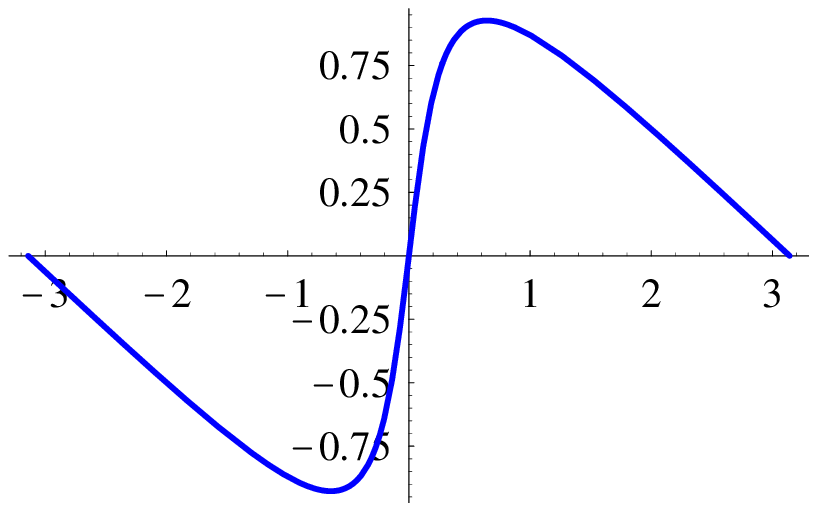}}
\end{center}
Furthermore, when $r$ varies from $0$ to $1$, $\iota(x)$ varies from
the constant zero function to a discontinuous linear function.  These
will be useful in calculations involving periodic orbits of $R_{r,\Omega}$.
For example, one may determine the birth of fixed point according to this
knowledge of $\iota$.
\begin{proposition}
For all $r$ and $\Omega \not\in [\pi-2a_r, \pi + 2a_r]$, the map
$R_{r,\Omega}$ has no fixed point and a saddle-node bifurcation occurs
at $\modulus{\Omega-\pi}=2a_r$, that is, when $r$ is the cosine of
the deviation of between incident and reflected angles.
\end{proposition}
\begin{proof}
To find fixed points of $R_{r,\Omega}$, one tries to solve the equation
$$
\tR_{r,\Omega}(x) = x + \Omega - 2\iota(x) = x \mod 2\pi,
$$
which can be reduced to
$$
\frac{\Omega}{2} = \iota(x) \mod \pi.
$$
Since $\iota(-a_r) \leq \iota(x) \leq \iota(a_r)$, for $\Omega \in
(-\pi,\pi]$, the equation has solution if and only if
$$
\frac{-\pi}{2}+a_r \leq \frac{\Omega}{2} \leq \frac{\pi}{2}-a_r .
$$
The fixed point of $\tR_{r,\Omega}(x)$ at the boundary parameter values
$\Omega = \pm(\pi-2a_r)$ occurs at $x = a_r$.  It is easy to see that
for $0 < r < 1$ and any $\Omega$, $\tR_{r,\Omega}'(a_r) = 1$ and
$\tR_{r,\Omega}''(a_r) \ne 0$, so as $(r,\Omega)$ crosses this boundary,
a saddle node bifurcation occurs.
\end{proof}
\begin{corollary}
Let $a_r$ as above and $b_r\in(0,a_r)$ such that $\cos b_r =
\dfrac{1+2r^2}{3r}$.  The region that $R_{r,\Omega}$ has
attracting fixed point is
$$
\Set{(r,\Omega) : 2\iota(b_r) < \modulus{\Omega} < 2\iota(a_r)}.
$$
In fact, the equation $2\iota(b_r) = \modulus{\Omega}$ determines
the generic values of $r$ for the happening of period-doubling
bifurcation. 
\end{corollary}
\begin{proof}
If $x\in(-\pi,\pi]$ corresponds to an attracting fixed point of
$R_{r,\Omega}$, we have $\tR_{r,\Omega}(x) = x \mod 2\pi$ and
$\modulus{\tR'_{r,\Omega}(x)} = \modulus{1-2\iota'(x)} < 1$.
From the expression of $\iota'(x)$, this is equivalent to
$$
\cos a_r = r < \cos x < \frac{1+2r^2}{3r}.
$$
Then, $x \in (-a_r,-b_r)\cup(b_r,a_r)$.  Since $\iota(x)$ is
increasing on both $(-a_r,-b_r)$ and $(b_r,a_r)$, $\Omega/2$
lies in the intervals defined by the image under $\iota$,
namely, $\iota(b_r) < \modulus{\Omega/2} < \iota(a_r)$.  To
see the period doubling bifurcation, it is sufficient to
show
$$
\left.
\frac{\partial(\tR_{r,\Omega}^2)'}{\partial r} \right|_{x=\pm
b_r} \ne 0, \qquad \text{along the curves $\Omega = \pm 2\iota(b_r)$}.
$$
Here we need some calculations which will also be useful in the future.
\begin{lemma}\label{lem-pdiffR}
We have
\begin{align*}
\frac{\partial \tR_{r,\Omega}(x)}{\partial r} &=
\frac{-2\sin x}{1 - 2r\cos x+r^2}, \\
\frac{\partial\tR_{r,\Omega}'(x)}{\partial r} &=
\frac{4r - 2(1+r^2)\cos x}{(1-2r\cos x+r^2)^2}.
\end{align*}
Furthermore,
$$
\frac{\partial(\tR^2_{r,\Omega})'(x)}{\partial r} =
\frac{\partial\tR_{r,\Omega}'(x)}{\partial r} \cdot
\tR_{r,\Omega}'(\tR_{r,\Omega}(x)) + \tR_{r,\Omega}'(x) \cdot
\left.\frac{\partial\tR_{r,\Omega}'}{\partial r}
\right|_{\tR_{r,\Omega}(x)} \cdot \frac{\partial
\tR_{r,\Omega}(x)}{\partial r}.
$$
In particular, if $\tR_{r,\Omega}(x)=x$ or
$\tR_{r,\Omega}(x)=-x$, one has
$$
\frac{\partial(\tR^2_{r,\Omega})'(x)}{\partial r} =
\tR_{r,\Omega}'(x) \cdot
\frac{\partial\tR_{r,\Omega}'(x)}{\partial r} \cdot \left[ 1 +
\frac{\partial \tR_{r,\Omega}(x)}{\partial r} \right].
$$
\end{lemma}
\begin{proof}[Proof of the lemma]
The first two results only require simple calculus.  The third
one is a repeated application of the chain rule.  Then,
using that both $\tR_{r,\Omega}'$ and
$\dfrac{\partial\tR_{r,\Omega}'}{\partial r}$ are even functions
in $x$, the last result follows.
\end{proof}
At $\Omega=\pm\iota(b_r)$, we have
$\tR_{r,\Omega}(b_r)=b_r$ and
$\tR_{r,\Omega}'(b_r)=-1$, therefore
$$
\left.
\frac{\partial(\tR_{r,\Omega}^2)'}{\partial r} \right|_{x=\pm
b_r} = - \left. \frac{\partial(\tR_{r,\Omega})'}{\partial r}
\right|_{x=b_r} \cdot \left[ 1 + \left.\frac{\partial
\tR_{r,\Omega}}{\partial r}\right|_{x=b_r} \right]
$$
Moreover, for $r>1/2$,
\begin{align*}
\left. \frac{\partial\tR_{r,\Omega}'}{\partial r}\right|_{x=b_r} &=
\frac{4r - 2(1+r^2)\cos b_r}{(1-2r\cos b_r+r^2)^2}
= \frac{-6(1-2r^2)}{r(1-r^2)} \ne 0. \\
\left. \frac{\partial \tR_{r,\Omega}}{\partial r}\right|_{x=b_r} &=
\frac{-2\sin b_r}{1 - 2r\cos b_r+r^2}
= \frac{-2\sqrt{4r^2-1}}{r\sqrt{1-r^2}} \ne -1, 
\end{align*}
except at $r = \sqrt{\frac{-15+\sqrt{241}}{2}} \approx 0.5119$.
Thus, generically, period-doubling bifurcation occurs at
$\Omega=\pm2\iota(b_r)$.
\end{proof}
The above information about fixed points is illustrated by
the following picture, in which the grey areas
are the region in $(r,\Omega)$-plane where fixed points exist.
The darker area is where attracting fixed points occur.
\begin{center}
\mbox{\epsfysize=5cm \epsfbox{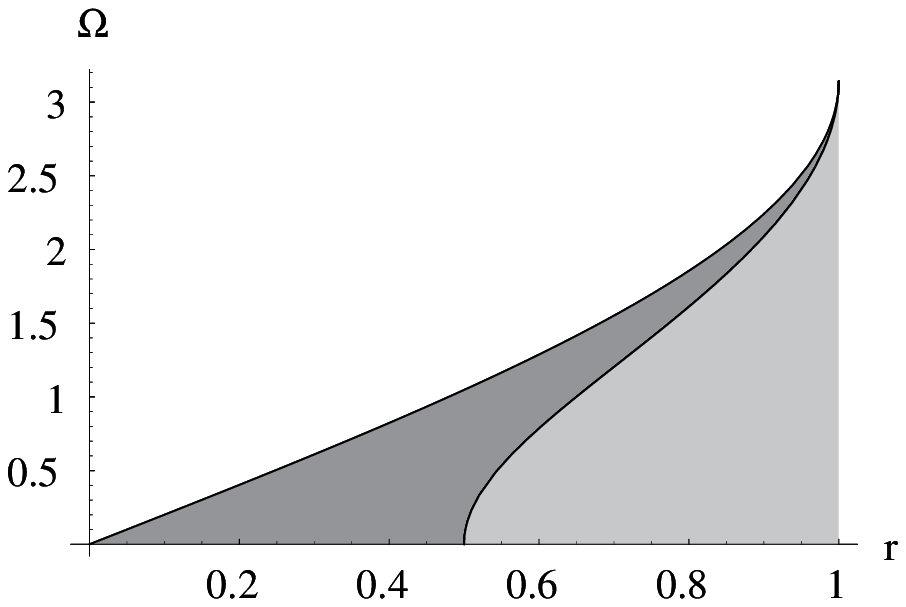}}
\end{center}

We would like to look into maps in the family that respect
the symmetry of the circle.  Let $\map{\rho}{\sphere^1}{\sphere^1}$
be the reflection across the real axis.  Then $R_{r,\Omega}\circ\rho
= \rho\circ R_{r,\Omega}$ if and only if $\Omega = 0, \pi$.
This symmetry corresponds to the fact that $\tR_{r,\Omega}(x)-\Omega$
is an odd function.

For $\Omega = 0, \pi$, if $n$ is the smallest integer that
$R^n_{r,\Omega}(\phi) = \rho(\phi) \ne
\phi$, one can easily show that $\phi$ belongs to a periodic orbit
of prime period~$2n$.  We call it a symmetric orbit of period $2n$.
In terms of $\tR_{r,\Omega}$, this corresponds to
$\tR_{r,\Omega}^n(x)=-x \ne x \mod 2\pi$.
An orbit is asymmetric of period $2n$ if $R^{2n}_{r,\Omega}(\phi)=\phi$
but $R^n_{r,\Omega}(\phi) \ne \rho(\phi)$.  If an asymmetric orbit
is formed by $\phi$, then another asymmetric one, called the twin orbit,
is formed by $\rho(\phi)$.  The twin orbit may be itself when
$\rho(\phi)=\phi$.  For the off-center reflection, a self-twin
asymmetric orbit must be the 2-cycle 
$\left\{e^{\thei 0}, e^{\thei\pi}\right\}$.
There are numerous articles about symmetric periodic orbits in
dynamical systems, especially on continuous types.  An early one is
\cite{D}.  Their attention is on the return map of some reversible
mechanic system such as the three-body system.

\begin{proposition}
For any $r>1/3$ and $\Omega = 0, \pi$, the map $R_{r,\Omega}$
has either no attracting orbit, or one symmetric attracting orbit,
or two (counting multiplicity) asymmetric twin attracting orbits.
\end{proposition}
\begin{proof}[Sketch of proof]
We mainly use two properties of $R_{r,\Omega}$ to conclude this. 
First, each $R_{r,\Omega}$ has negative Schwarzian derivative,
the technique of \cite{CE} or \cite{Br} is applicable.  Therefore,
once the map $\tR_{r,\Omega}$ has an attracting cycle, at least
one critical point of it must be attracted to this attracting
orbit.  Since $R_{r,\Omega}$ has only two critical points, there are
at most two attracting orbits.  If the attracting orbit is a
reflection symmetric one, then this orbit attracts both critical points. 
If the orbit is an asymmetric one, its twin orbit attracts another
critical point.  Multiplicity occurs if there is a self-twin orbit.
\end{proof}
Besides the above general information about symmetric and
asymmetric cycles, it is interesting to know two specific cases
about it.  The first is about self-twin asymmetric 2-cycles.
\begin{proposition}
A single asymmetric 2-cycle for $R_{r,\Omega}$ occurs if and
only if $\Omega=\pi$.  The orbit is $\Set{e^{\thei 0},
e^{\thei\pi}}$ which is attracting when $r<1/\sqrt{5}$.
\end{proposition}
\begin{proof}
The orbit for a single asymmetric 2-cycle can be obtained by
solving the equation
$$
\frac{\Omega\pm\pi}{2} = \iota(x) \mod \pi .
$$
The only simultaneous solution exists when $\Omega=\pi$ and
$\phi=e^{\thei 0}, e^{\thei\pi}$.  Then, whether the orbit
is attracting can be easily decided by calculating
$\tR'_{r,\pi}(0)\tR'_{r,\pi}(\pi)$.
\end{proof}
Secondly, one expects that symmetric cycles occur naturally
for $\Omega = 0, \pi$.  Do such cycles exist even if the map
is not ``symmetric''?  The following result provides a partial
evidence for the answer.
\begin{proposition}
The map $R_{r,\Omega}$ has a symmetric 2-cycle if and only if
$\Omega = 0, \pi$.  In both cases, the symmetric 2-cycle is
unique.  The cycle of $R_{r,0}$ is
attracting while that of $R_{r,\pi}$ is repelling.
\end{proposition}
\begin{proof}
Let $\Set{e^{\thei x}, e^{-\thei x}}$ be a symmetric 2-cycle, i.e.,
$R_{r,\Omega}(e^{\pm\thei x}) = e^{\mp\thei x}$.  Equivalently,
\begin{align*}
x + \frac{\Omega}{2} &= \iota(x) \mod \pi, \\
\intertext{and also,}
x - \frac{\Omega}{2} &= \iota(x) \mod \pi.
\end{align*}
If we subtract the second equation from the first one,
we obtained that $\Omega = 0 \mod \pi$ and $\Omega \in (-\pi,\pi]$.
Clearly, $\Omega = 0, \pi$.

For $\Omega = 0$, the trivial solutions $x = 0, \pi$ corresponds
to fixed points of $R_{r,0}$.  Further calculation on $\iota$
shows that the equation~$x = \iota(x)$ has only a nontrivial solution
for $r > 1/2$.  In fact, let $c_1 \in  (0,\pi/2)$ such that
$\cos c_1 = \dfrac{1}{2r}$, then $\pm c_1$ form a symmetric 2-cycle.
Since $\tR_{r,0}'(c_1)\tR_{r,0}'(-c_1) =
\left[\dfrac{1-3r^2}{r^2}\right]^2$, it follows that the cycle
is attracting for $\dfrac{1}{2} < r < \dfrac{1}{\sqrt{2}}$.
Moreover, since $\modulus{\iota(x)} \leq \pi/2$, one can only
obtain trivial solution $0, \pi$ from $x = \iota(x) + k\pi$ for
$k \ne 0$.

For $\Omega=\pi$, the above equations
always has solution for all $r$.  In fact, its solution within
$(-\pi,\pi)$ is given by the following.
Let $c_2 \in (\pi/2,\pi)$ satisfy $\cos c_2 =
\dfrac{1-\sqrt{1+8r^2}}{4r}$.  We claim that $c_2$ and $2\pi-c_2$
are solutions to the equation~$x - \frac{\pi}{2} = \iota(x)$ and
$x + \frac{\pi}{2} = \iota(x)$ respectively.  Taking tangent to both
sides and using that $\iota(x) = \Arg(\cos x-r+\thei\sin x) -
\Arg(\cos x + \thei\sin x)$, we have
$$
\frac{-\cos x}{\sin x} = \frac{r\sin x}{1-r\cos x}.
$$
It follows that $2r\cos^2 x - \cos x - r = 0$.  Hence
$\Set{e^{\thei c_2},e^{-\thei c_2}}$ is a symmetric 2-cycle.
It can be easily checked that
$$
\tR_{r,\pi}'(c_2) = \tR_{r,\pi}'(2\pi-c_2) =
\frac{2(\sqrt{1+8r^2}+3r^2)}{1 + \sqrt{1+8r^2} + 2r^2} > 1.
$$
Therefore, the symmetric 2-cycle is repelling.  Again, by that
$\modulus{\iota(x)}\leq \pi/2$, we do not have other solutions
to $x \pm \frac{\pi}{2} = \iota(x) + k\pi$ for $k \ne 0$.
\end{proof}

\section{Bifurcation}
In this last section, our aim is to understand bifurcations
between symmetric and asymmetric orbits of this family of maps.
In particular, we would like to address some of the questions
on the dynamics of $R_{r,\pi}$ raised in \cite{AL}.  Let us first
look at the asymptotic orbit diagrams of the critical points of
$R_{r,0}$ and $R_{r,\pi}$, in which the bifurcations are shown.

\noindent
\begin{center}
%\vspace*{5cm}
\mbox{\epsfysize=4.7cm \epsfbox{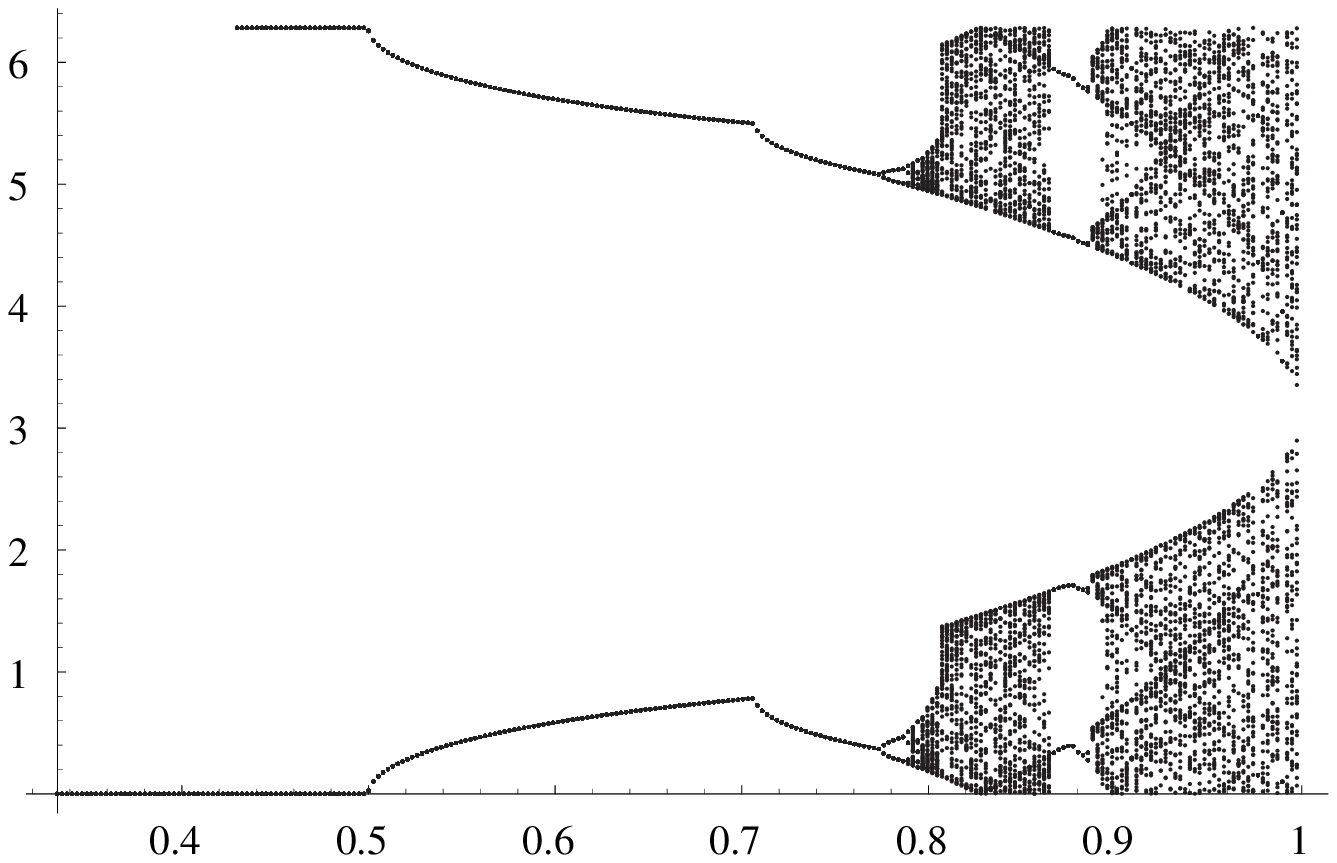}}\hfil
\mbox{\epsfysize=4.7cm \epsfbox{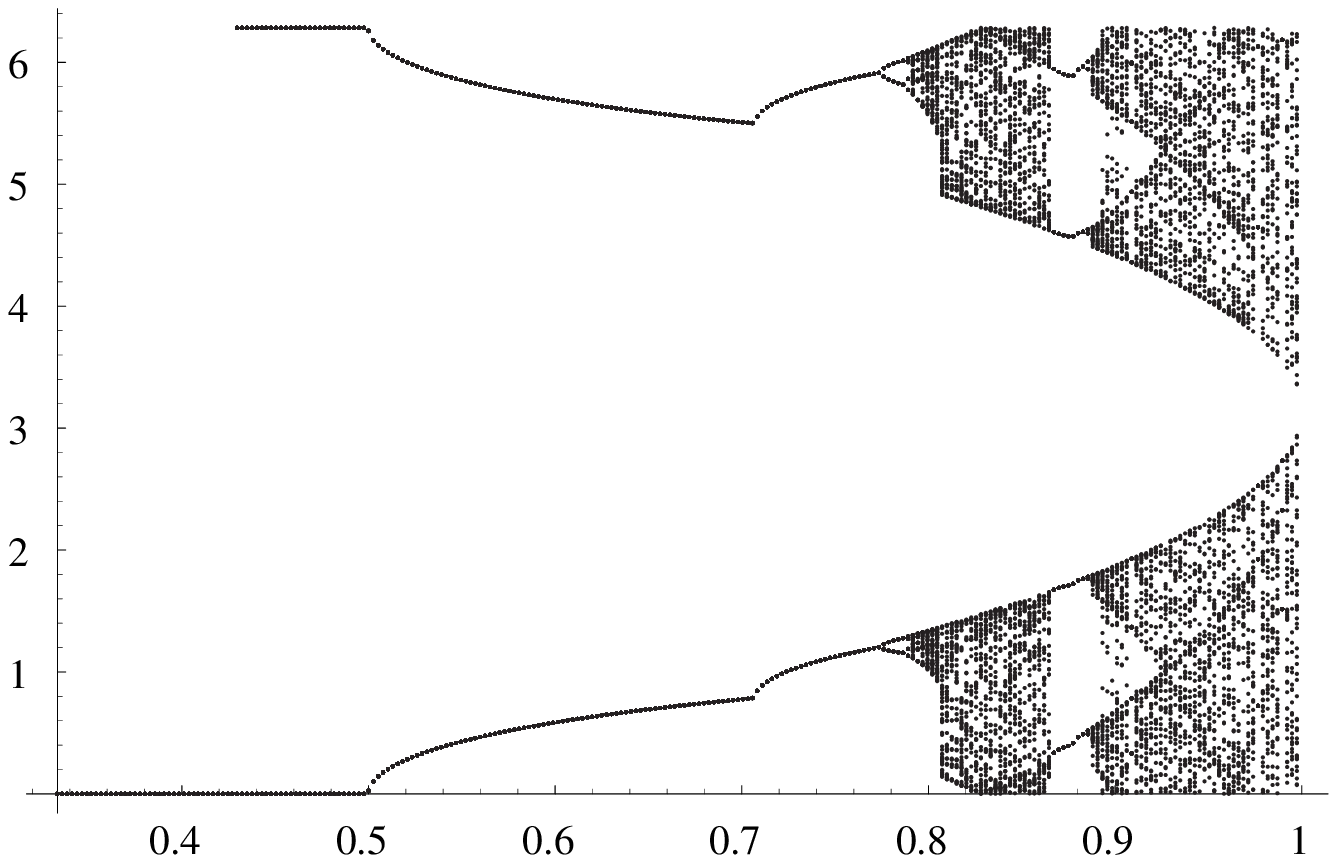}} \\
{\footnotesize Asymptotic orbits of critical points of $R_{r,0}$
(low resolution).}
\end{center}

\bigskip

\noindent
\begin{center}
%\vspace*{5cm}
\mbox{\epsfysize=5cm \epsfbox{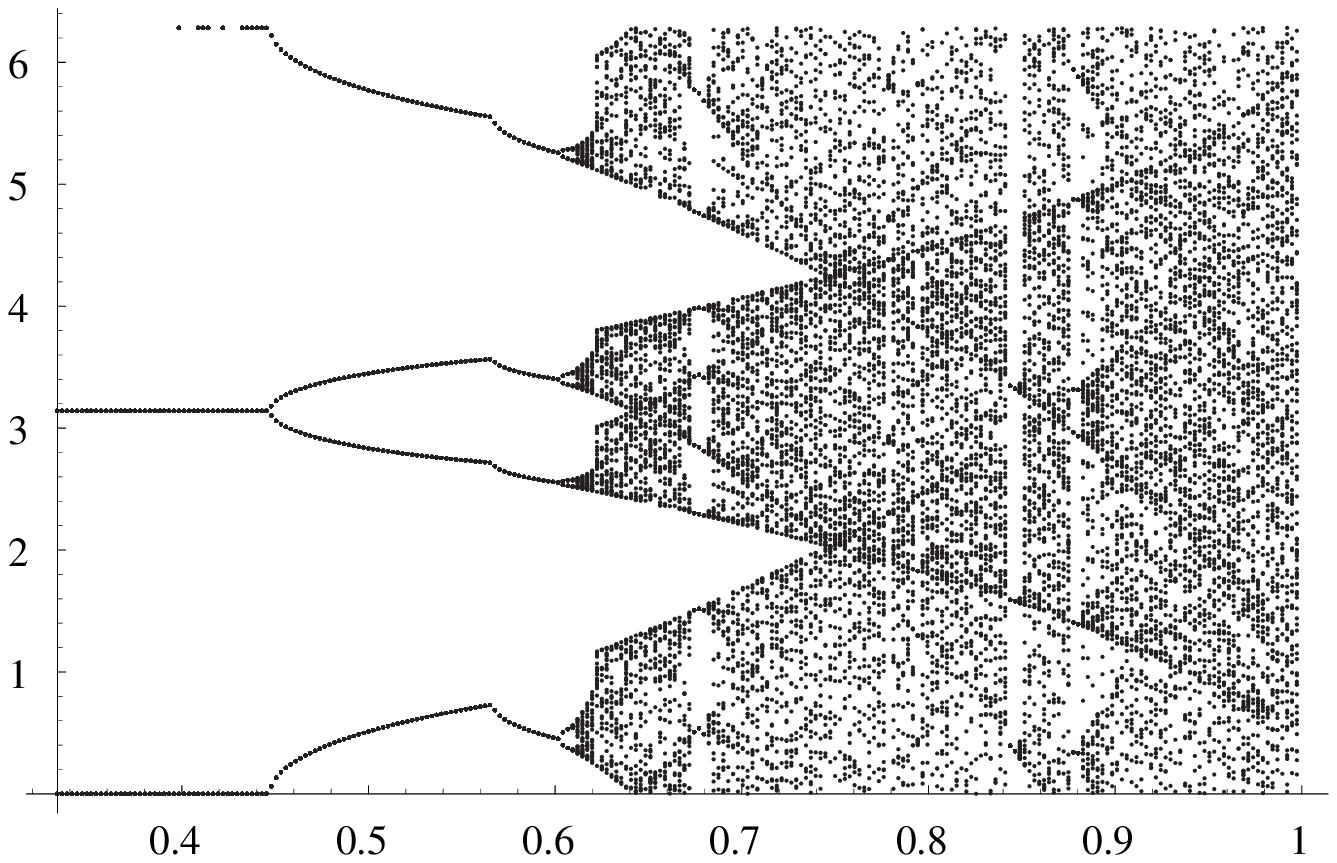}}\hfil
\mbox{\epsfysize=5cm \epsfbox{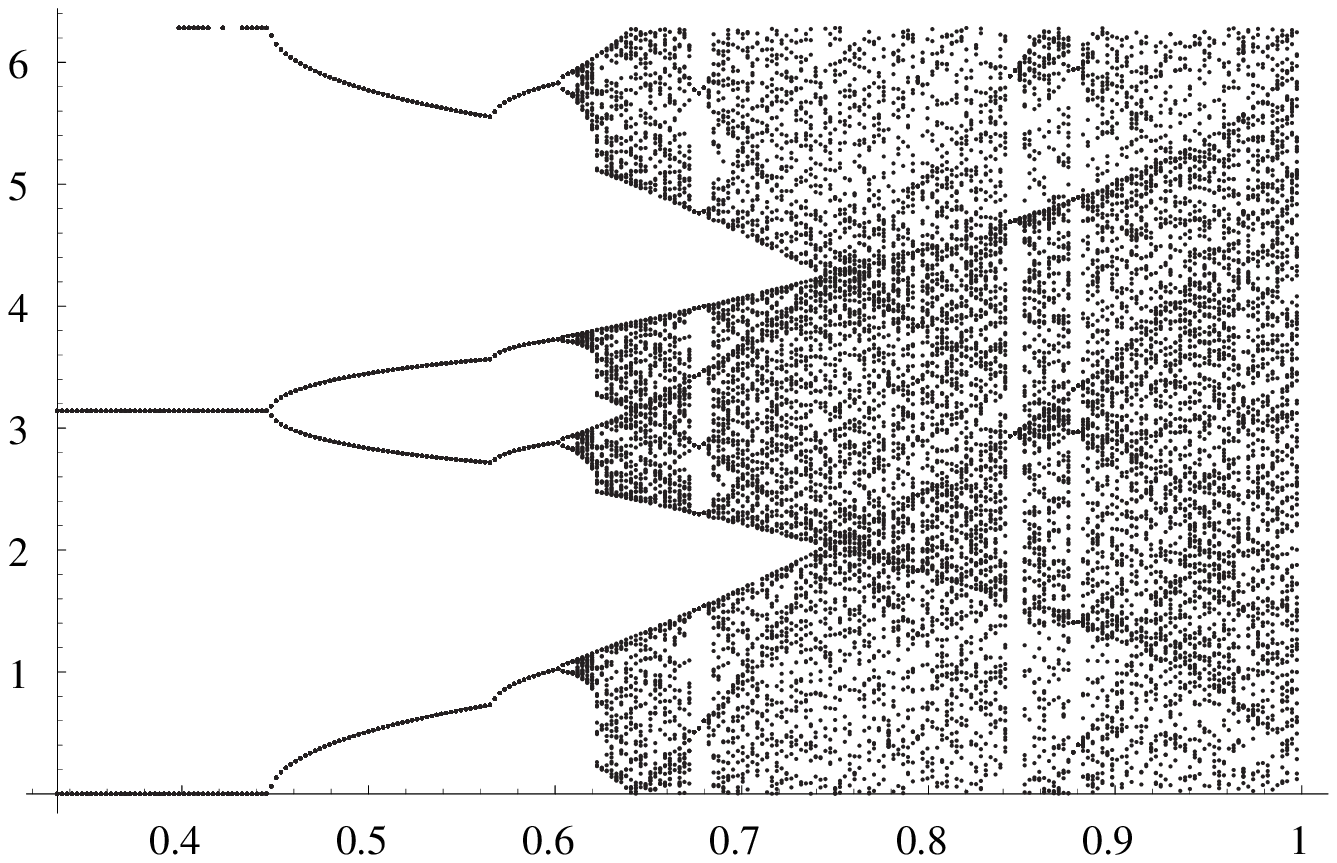}} \\
{\footnotesize Asymptotic orbits of critical points of $R_{r,\pi}$
(low resolution).}
\end{center}
In both pictures, there is an obvious complementary nature and
there are half-branches of bifurcations.  These will be explained
analytically in the coming propositions.
\begin{proposition}
The asymmetric 2-cycle of $R_{r,\pi}$ bifurcates into a symmetric
attracting 4-cycle at $r=1/\sqrt{5}$.  For $R_{r,0}$, there is
a pitch-fork period preserving bifurcation of the symmetric
2-cycle into asymmetric attracting ones at $r = 1/\sqrt{2}$.
\end{proposition}
\begin{proof}
For the 2-cycle $\Set{e^{\thei 0}, e^{\thei\pi}}$ of $R_{r,\pi}$,
the derivatives are given by
\begin{align*}
\tR'_{r,\pi}(0) &= \frac{1-3r}{1-r} ;\\
\tR'_{r,\pi}(\pi) &= \frac{1+3r}{1+r}.
\end{align*}
Thus, the 2-cycle undergoes a period-doubling bifurcation when
$-1 = \dfrac{1-9r^2}{1-r^2}$, which is exactly $r=1/\sqrt{5}$.
We claim that the new attracting 4-cycle is reflection symmetric.
In fact, the equation for a reflection symmetric 4-cycle of
$R_{r,\pi}$ is
$$
x + \pi = \iota(x) + \iota\circ\tR_{r,\pi}(x) \mod \pi.
$$
Let $f(x) = x - \iota(x) - \iota(x+\pi-2\iota(x)) + \pi$.  Then
we have to solve for $f(x) = 0 \mod \pi$.  Clearly, $f(x)-\pi$
is an odd function with $f(-\pi) = 0$, $f(0) = \pi$, and $f(\pi)
= 2\pi$.   If $f$ is decreasing at $-\pi$, equivalently, it is
so at $\pi$, then $f$ must have a zero modulo~$2\pi$
in neighborhoods of $\pm\pi$.  It is easy to compute that
\begin{align*}
f'(-\pi) = f'(\pi) &= 1 - \iota'(-\pi) -\iota'(0)
\left[ 1-2\iota'(-\pi)\right] \\
&= \frac{1-5r^2}{1-r^2}.
\end{align*}
Thus, $f'(-\pi) \leq 0$ if and only if $r\geq 1/\sqrt{5}$.  This
shows that for $1/\sqrt{5} < r$, $R_{r,\pi}$ has a symmetric
4-cycle.  Moreover, it must be attracting when $r <
1/\sqrt{5}+\varepsilon$ by continuity.

We have shown in the previous section that an attracting
symmetric 2-cycle exists for $R_{r,0}$ with $1/2 < r < 1/\sqrt{2}$.
It is given by
$
\tR_{r,0}(c_1) = -c_1 \mod \pi.
$
One may further calculate according to the lemma in \S 2
to obtain that, at $r=1/\sqrt{2}$ and $x=c_1$,
$$
\frac{\partial^2 R^2_{r,0}}{\partial r \partial x} =
2\sqrt{2}-8 \ne 0.
$$
Then, by an argument making use of the Inverse Function Theorem,
one concludes that $R_{r,0}$ has a 2-period preserving pitch
fork bifurcation at $r=1/\sqrt{2}$.
\end{proof}

\begin{proposition}\label{symbreak}
There is a pitch-fork bifurcation on $R_{r,\pi}$ where a symmetric
orbit of period 4 breaks into two asymmetric orbits of period 4.
\end{proposition}
\begin{proof}[Sketch of proof]
The key is to consider the zero set of $\tR^4_{r,\pi}(x) = x$ in the
$(r,x)$-plane.  In the above, we have shown that for $r$ slightly larger
than $1/\sqrt{5}$, this can be obtained from the solution set of
a symmetric 4-cycle, $\tR^2_{r,\pi}(x)=-x$.  We then solve for the
specific values of $(r_0,x_0)$ such that $\tR^2_{r_0,\pi}(x_0)=-x_0$
and $(\tR^2_{r_0,\pi})'(x_0) = 1$.  
Our calculation involves the expression of $R_{r,\pi}$ as a Blaschke
product given in \cite{AL},
$$
R_{r,\pi}(z) = \frac{-z^2(1-rz)}{z-r}.
$$
Then, a symmetric 4-cycle can be solved from $z=e^{\thei\phi}$ and
$$
R_{t,\pi}^2(z) = z^4\cdot \frac{(1-rz)^2}{(z-r)^2} \cdot
\frac{r^2z^3-rz^2-z+r}{rz^3-z^2-rz+r^2} = \frac{1}{z}.
$$
After factoring out the obvious factors $(z+1)(z-1)$, and letting 
$y=\cos\phi$, we have the polynomial equation
$$
-1-4r^2+r^4 + 2r(1+7r^3)y + 4r^2(2-3r^2)y^2 - 24r^3y^3 + 16r^4y^4 =
0.
$$
For $r>1/\sqrt{5}$, all the four roots of this equation are real.
Two of them are always~$> 1$, one lies within $(-1,0)$ and the other
within $(0,1)$.  After taking arccosines, the four solutions form
a symmetric 4-cycle.  Let $x_0=x_0(r)$ be a solution, then we
solve for $r$ in $(\tR^2_{r,\pi})'(x_0) = 1$,
we obtain a numerical value of
$r_0$ approximately equal to 0.57.  We further verify that
$$
\pdiff{(\tR^4_{r,\pi})'}{r}(x_0,r_0) \ne 0.
$$
Many of the above calculations are lengthy and we indeed make use of
computation software to help us.  Finally, we apply
the Inverse Function Theorem to conclude that there is a 4-period
preserving pitch fork bifurcation.
\end{proof}
The local change of the graphs of $R^4_{r,\Omega}$ is shown in
the picture.

\noindent
\begin{center}
%\vspace*{5cm}
\mbox{\epsfysize=5cm \epsfbox{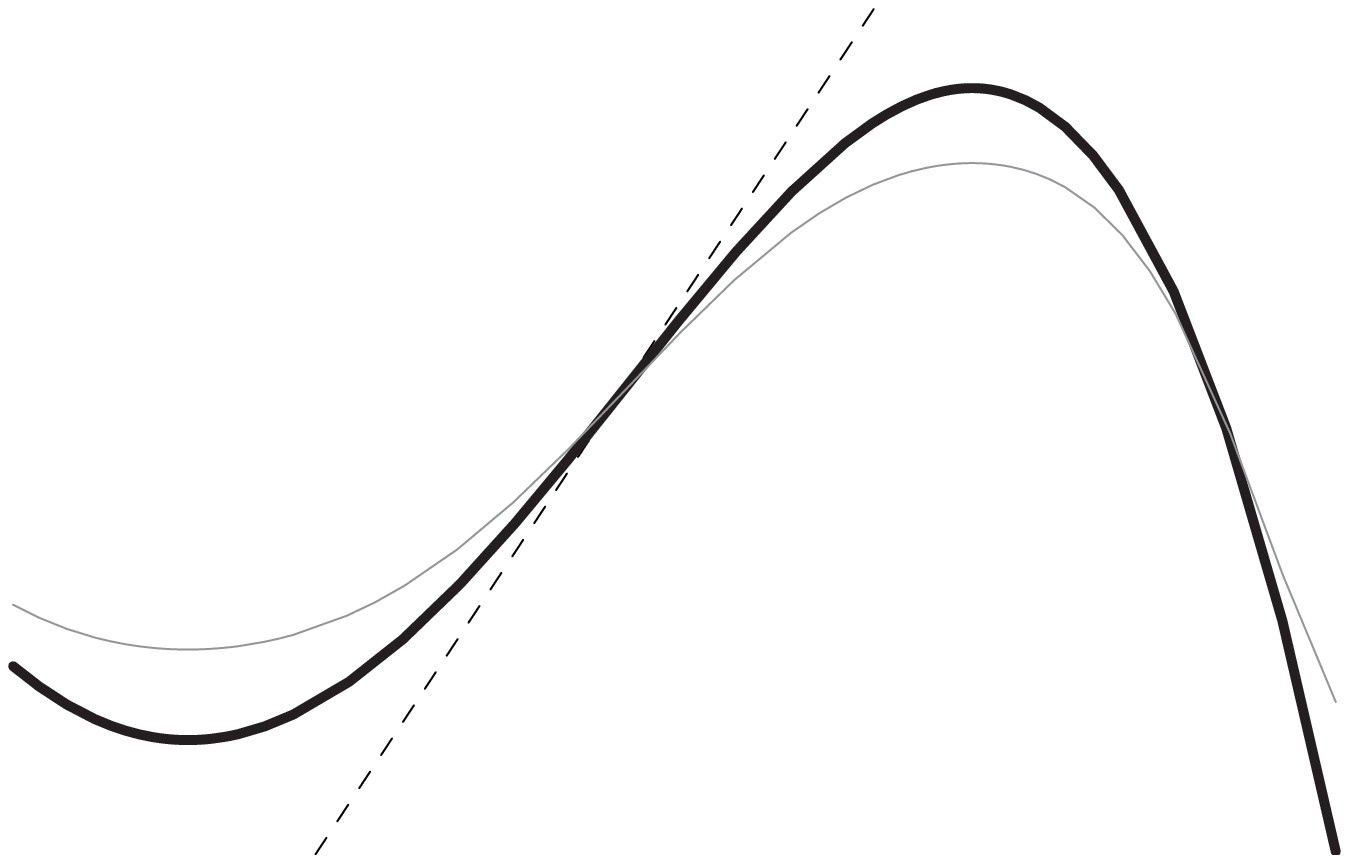}}\hfil
\mbox{\epsfysize=5cm \epsfbox{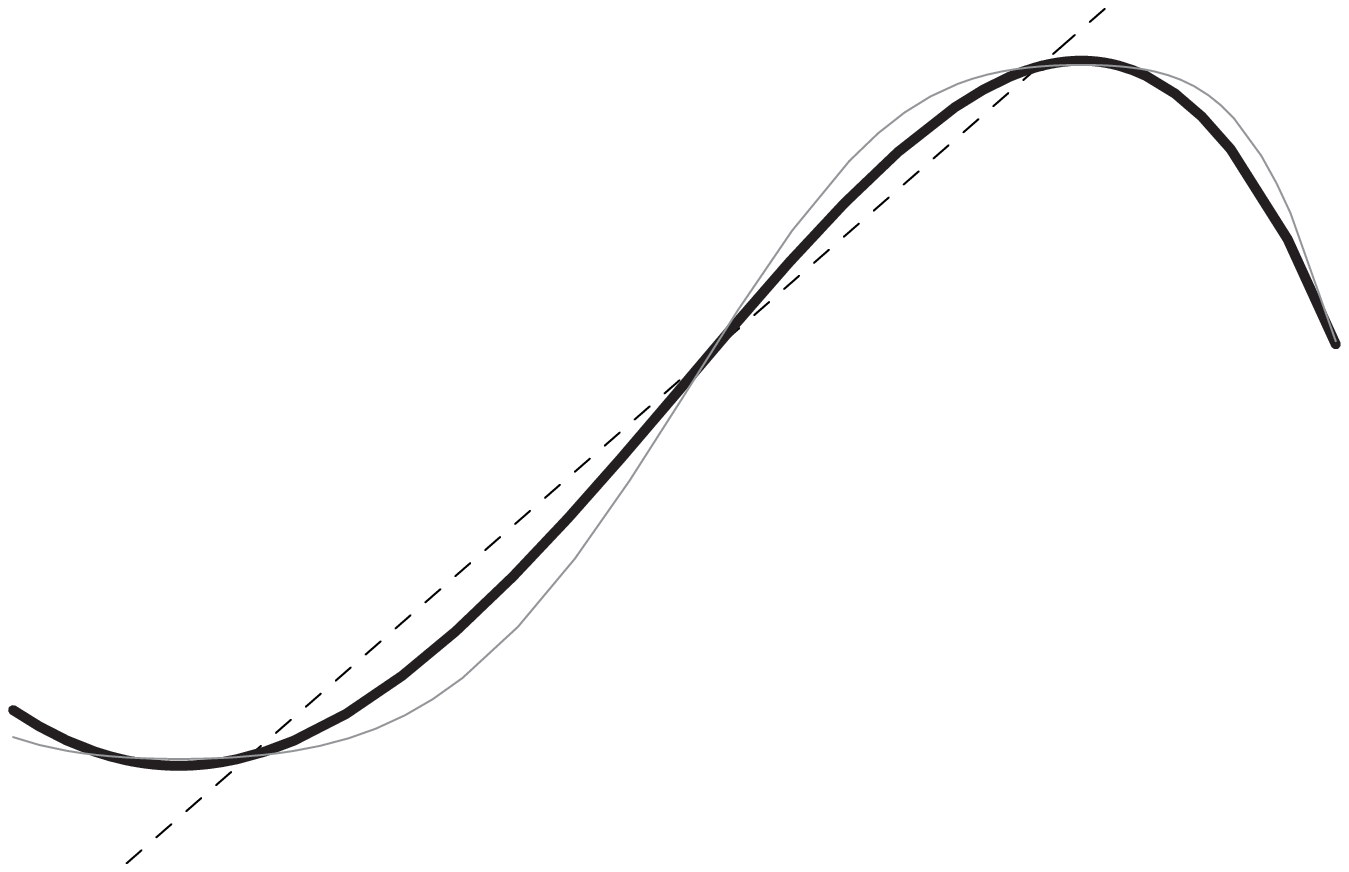}} \\
{\footnotesize Graphs of $R^4_{r,\Omega}$ (darker) and $R^8_{r,\Omega}$ (lighter) before and after the pitch-fork bifurcation.}
\end{center}

To end this article, we would like to present two pictures,
asymptotic orbits for $\Omega = \pi/2, \pi/4$ to illustrate
the wide variation of the dynamics of this family with respect
to $\Omega$.

\noindent
\begin{center}
%\vspace*{5cm}
\mbox{\epsfysize=5cm \epsfbox{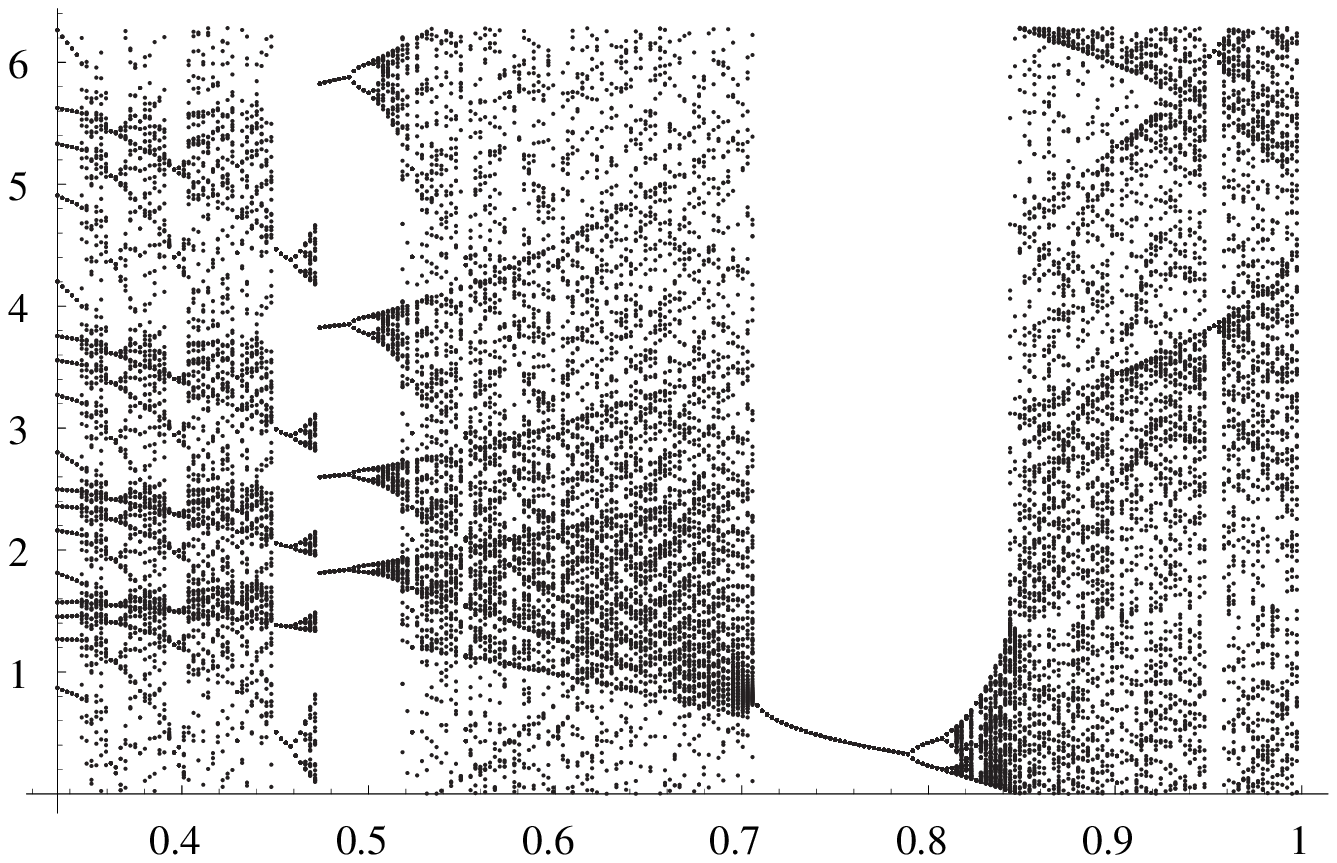}}\hfil
\mbox{\epsfysize=5cm \epsfbox{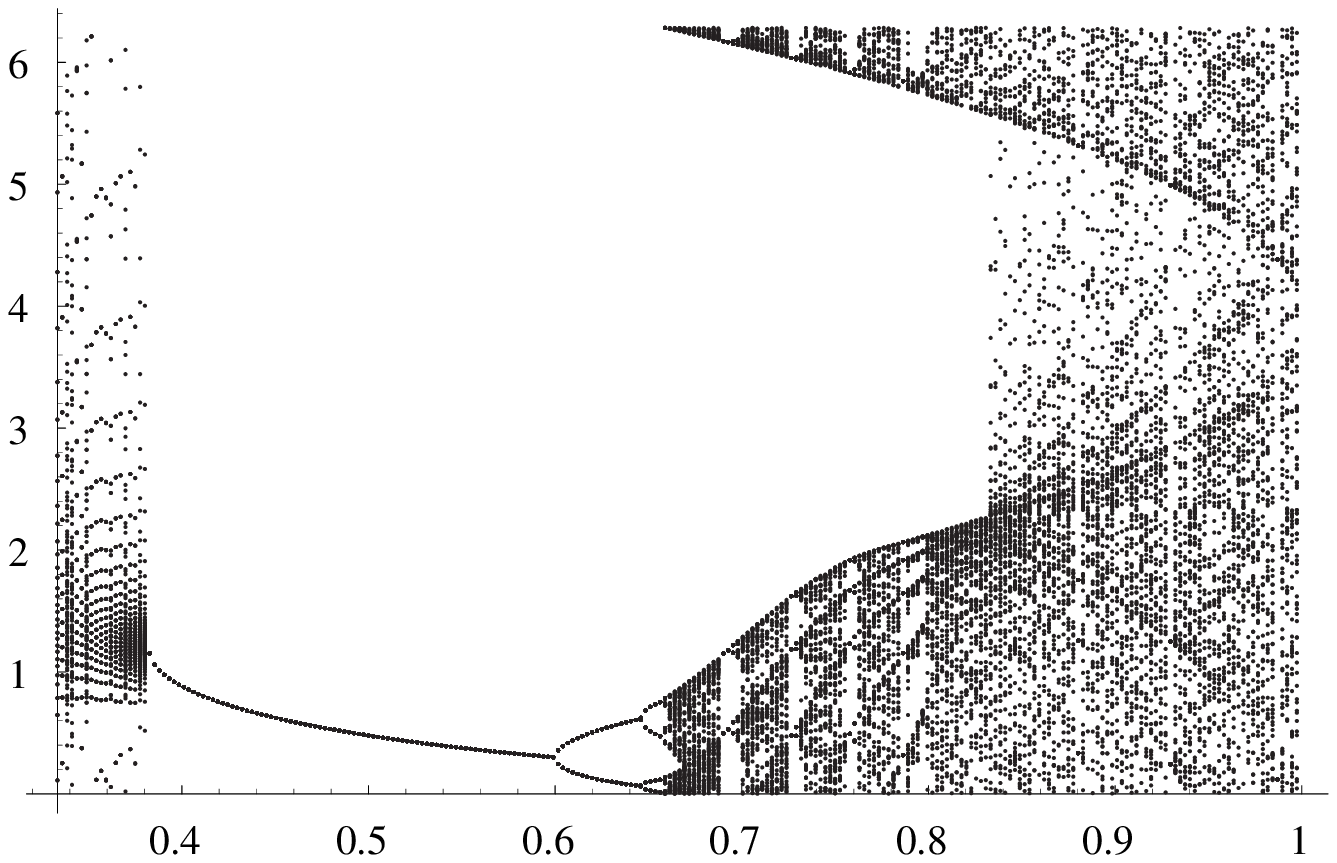}} \\
{\footnotesize Asymptotic orbits of $R_{r,\pi/2}$ and $R_{r,\pi/4}$
(low resolution).}
\end{center}

\end{document}